\documentclass[12pt]{article}

\usepackage{amssymb}
\usepackage{amsfonts}
\usepackage{latexsym}

\newcommand{\N}{{\mathfrak n}_f}
\newcommand{\DC}{{\mathcal C}_{\{M_j\}}(J)}

\newcommand{\DCA}{{\mathcal C}_{A}([0,1])}

\newcommand{\R}{\mathbb R}

\title{Lower bounds for quasianalytic functions, {\rm I}.
\\ How to control smooth functions?}
\date{}

\author{F. Nazarov, M. Sodin\thanks
{Supported by the Israel Science Foundation of the Israel Academy
of Sciences and Humanities under Grant No. 37/00-1.}, \ A. Volberg}

\begin{document}
\maketitle

\centerline{\em In memory of Th$\mbox \o$ger Bang (1917--1997)}

\begin{abstract}
Let $\mathcal F$ be a class of functions with the uniqueness
property: if $f\in \mathcal F$ vanishes on a set $E$ of positive
measure, then $f$ is the zero function. In many instances, we
would like to have a quantitative version of this property, e.g. the
estimate from below for the norm of the restriction operator
$f\mapsto f\big|_E$ or, equivalently, a lower bound for $|f|$ outside a
small exceptional set.
Such estimates are well-known and useful for polynomials,
complex- and real-analytic functions, exponential polynomials. In
this work we prove similar results for the Denjoy-Carleman and the
Bernstein classes of quasianalytic functions.

In the first part, we consider quasianalytically smooth functions.
This part relies upon Bang's approach and includes the proofs of
relevant results of Bang.  In the second part, which is to be
published separately, we deal with classes of functions
characterized by exponentially fast approximation by polynomials
whose degrees belong to a given very lacunar sequence.

The proofs are based on the elementary calculus technique.
\end{abstract}

\bigskip
\centerline{\sc \S1. Motivation and the results}

\bigskip Let $P$ be a polynomial. Its degree $d$ governs
the behaviour of $P$ on any interval $I\subset \R$, for instance,
$P$ {\em has at most $d$ zeroes on $I$} and
{\em for any measurable subset $E\subset I$}
$$
||P||_I \le \left( \frac{4|I|}{|E|}\right)^d
||P||_E\,.
\eqno (1.1)
$$
Here and everywhere below, $||f||_K=\sup_K|f|$ denotes the
supremum norm on $K$, and $ |E| $ denotes the length of a set $ E\subset 
\mathbb R $.

The first fact hardly needs any comments. The second one is a
rough version of the classical Remez inequality \cite{Remez} (cf.
\cite{DR} and \cite{BG}). Different proofs of (1.1) are known. The
simplest one uses the Lagrange interpolation formula for $P$ with
$d +1$ nodes on $E$ spaced by at least $|E|/d$, though this
gives us (1.1) with a worse constant $2e$ instead of $4$ on the
right hand side \cite{Nadirashvili}. The Remez inequality has a
plenty of applications and extensions, some of them can be found
in \cite{DR}, \cite{BG}, \cite{Nadirashvili},
\cite{Nadirashvili1}; by no means is this list complete.
The inequality is sufficiently sharp to capture that $P$ cannot have
zeroes of multiplicity larger than $d$.

Turning to analytic functions, we encounter another quantity which
controls their behaviour. Let $G\subset \mathbb C^1$ be a bounded
domain, $K\subset G$ be a compact subset, and let $f$ be a
bounded analytic function in $G$. Then the logarithm of the ratio
$$
\mathfrak B_f(K, G) = \log \frac{||f||_G}{||f||_K}
$$
is called the {\em Bernstein degree} of $f$ on $(K, G)$. If $P$ is
polynomial of degree $d$, and $G_R\subset \mathbb C^1$ is the
ellipse with the foci at $-1, +1$ and the semiaxes $R$, then by
the classical Bernstein inequality
$$
\mathfrak B_P([-1,1], G_R) \le d \log R.
$$

The Bernstein degree controls the number of zeroes of $f$ on $K$
as well as the local oscillations of $f$. By the Jensen formula
the number of zeroes of $f$ on $K$ counting with multiplicities
does not exceed $\gamma (K, G) \mathfrak B_f(K, G)$ where $\gamma
(K, G)$ depends only on the geometry of the couple $(K, G)$.
The Cartan lemma yields local estimates on $K$ similar
to (1.1) with the exponent $\gamma (K, G) \mathfrak B_f(K, G)$.
The interest to this classical theme was recently revived (cf.
\cite{Brudnyi}, \cite{IYa}, \cite{Il}, \cite{RY}, \cite{Y} and the
references therein).

In this work we shall exhibit a new index which controls
in a similar fashion the behaviour of quasianalytically smooth
functions.

Given closed interval $J\subset \R$ and given a sequence of
positive numbers $\{M_j\}$, introduce the class $\DC$ of
$C^\infty(J)$-functions such that
$$
||f^{(j)}||_J \le M_j, \qquad j\in \mathbb Z_+.
\eqno (1.2)
$$
We assume that the sequence $\{M_j\}$ is {\em logarithmically
convex}, that is
$$
M_j^2 \le M_{j-1} M_{j+1}, \qquad j\in \mathbb N. \eqno
(1.3)
$$
A convenient way to generate logarithmically convex sequences is
to fix a non-decreasing function $A: [1,\infty)\to (0,\infty)$ and
set $$M_j= M_{j-1}A(j), \qquad j\ge 1. \eqno (1.4)$$ Rescaling the
argument of $f$ and multiplying $f$ by a constant, we can always
assume that $f$ is defined on the interval $J=[0,1]$ and that
$M_0=1$. Usually, we shall keep this normalization and denote the
normalized Denjoy-Carleman classes by $\DCA$.

According to the classical Denjoy-Carleman theorem \cite{Carleman}
divergence of the series
$$
\sum_{j=1}^\infty \frac{M_{j-1}}{M_j} = +\infty \eqno (1.5)
$$
(or equivalently of the integral $$\int_1^\infty
\frac{dt}{A(t)}=+\infty \Big) $$ is a necessary and sufficient
condition for
quasianalyticity of the class $\DC$ \cite{Carleman} (that is,
$\DC$ contains no non-trivial function which vanishes at a point
with all derivatives). In the paper \cite{Bang} published 50 years
ago, Bang gave an intrinsic and elementary real variable proof of
the uniqueness part of the Denjoy-Carleman theorem. Strangely
enough, this concise paper left no trace in the vast literature
devoted to quasianalytic functions, unlike Bang's thesis
\cite{Bang0} which appeared to be more influential (cf.
\cite[Chapter~IV]{Mandelbrojt}, \cite{Cohen},
\cite[Section~1.3]{Hormander}). For this reason, we took a liberty
to reproduce (with minor variations) some results from \cite{Bang}
with their proofs.

\medskip\par\noindent{\bf Definition.} The {\em Bang degree} $\N$ of the
function $f\in \DCA$ is the largest integer $N$ such that
$$
\sum_{\log ||f||^{-1}_{[0,1]}<j\le N} \frac{M_{j-1}}{M_j} < e
\eqno (1.6)
$$
If the set of positive integers $N$ satisfying (1.6) is unbounded,
then we set formally $\N=+\infty$.

\medskip After a  minute reflection, one can see a certain similarity
between the Bernstein and the Bang degrees. The latter depends on
the growth of the sequence $\{M_{j-1}/M_j\}$ (that is, on the a
priori smoothness of $f$), and on the lower bound for
$||f||_{[0,1]}$ (the closer is $||f||_{[0,1]}$ to its a priori
upper bound $M_0=1$, the smaller is the degree $\N$). If the
series (1.5) diverges, that is the class $\DCA$ is quasianalytic,
then the degree $\N$ is always finite. In the non-quasianalytic
case, the degree can be infinite. In fact, we can allow the
function $f$ to have only finite smoothness: if $f\in C^m([0,1])$,
then we simply put $ A(j)=+\infty $ starting with $j=m+1$.

\medskip\par\noindent{\bf Theorem A }(Bang \cite{Bang}). {\em
The total number of zeroes (counting with multiplicities) of any
function $f\in \DCA$ does not exceed its Bang degree $\N$.
}

\medskip The first result of that type was
conjectured by Borel and proved by Carleman in \cite[p.24--27]{Carleman}.
Carleman based the proof on the Fourier transform and harmonic estimation.
The theorem says that if $f\in
C^n([0,1])$ and satisfies
$$
f^{(j)}(0)=0, \quad 0\le j\le n-1, \qquad \hbox{and} \quad f(1)=1,
$$
then
$$
\sum_{j=1}^n \frac1{M_j^{1/j}} < 8\left(\frac1{2} + \pi e + 2
\sqrt{\pi e} \right)
$$
where $ M_j = ||f^{(j)}||_{[0,1]}$.

This estimate cannot be deduced directly from Theorem~A since the
sequence $\{M_j\}$ we deal with is assumed to be logarithmically
convex whereas Borel and Carleman did not impose any condition on
that sequence. Nevertheless, as we shall see in subsection 5.4, there is a
more general version of Bang's result which contains the result of
Carleman.

One can probably extract from Hirschman's paper \cite{Hirschman} a result
similar to Theorem~A (even with $\frac{2}{\pi}$ instead of $e$ on the
right hand side of (1.6)\,), however with some additional regularity of
the sequence $\{M_j\}$. Hirschman used the Carleman technique combined
with the Cartan-Gorny estimates of derivatives of smooth functions.

In the second theorem it will be convenient to assume that the
function $A$ which defines according to (1.4) the sequence
$\{M_j\}$ is a $C^1$-function (if $A$ is a piecewise linear function,
then one can use the left derivative). We set
$$
\gamma (n) := \sup_{1\le s \le n} \frac{sA'(s)}{A(s)} \eqno (1.7)
$$
and
$$
\Gamma (n)=4e^{4+ \gamma (n)}.
$$

\medskip\par\noindent{\bf Theorem B. }{\em
Suppose $f\in \DCA$. Then for any interval $I\subset [0,1]$ and
any measurable subset $E\subset I$
$$
||f||_I \le \left( \frac{\Gamma (2\N)|I|}{|E|}\right)^{2\N}
||f||_E\,. \eqno (1.8)
$$
}

\medskip We say that the class $\DCA$ is {\em regular} if the
supremum
$$
\gamma := \sup_{s\ge 1} \frac{sA'(s)}{A(s)}
$$
is finite. For example, the real analytic class ($A(s)=s$) and the
logarithmic classes ($A(s)=s\log^\alpha (s+e)$) are regular. For regular
classes, estimate (1.8) holds with the factor $\Gamma =4e^{4+ \gamma}$ on
the right hand side.

Theorems A and B show that Bang's degree is an important
characteristics of smooth functions. However, we do not know much
about its basic properties. For example, if $ f $ is a polynomial,
how to bound from above the Bang degree $\N$ by the usual
degree? If $ f $ is real analytic, the same question can be asked
about the upper bound of the Bang degree by the Bernstein degree.
Recently, N. Roytvarf \cite{R} and D. Novikov and S. Yakovenko
\cite{NYa} obtained useful estimates for the Bernstein degree of
linear combinations, products, (analytic) quotients and
derivatives of given functions. It seems to be interesting to get
results in that spirit for the Bang degree. At last, it looks
probable, that Bang's degree has a certain invariance under real
analytic diffeomorphisms of the interval $ [0,1] $.

\medskip\par\noindent{\sc Acknowledgment. } The authors thank
Alexander Borichev for numerous useful remarks.

\bigskip\centerline{\sc \S 2. Bang's fundamental inequality}

\bigskip
Given a function $f\in \DCA$ and a point $x\in [0,1]$, we define
the {\em norm of $f$ at $x$} as
$$
B_f(x):= \max_{j\ge 0} \frac{|f^{(j)}(x)|}{e^j M_j}.
$$
A small norm means that a large section of the sequence
$\left\{\frac{|f^{(j)}(x)|}{M_j}\right\}_{j\ge 0}$ consists of
small numbers. For example, $B_f(x) \le e^{-q}$ for some $q\in
\mathbb Z_+$ if and only if $|f^{(j)}(x)|\le e^{j-q}M_j$ for $0\le
j\le q$.

\medskip\par\noindent{\bf Lemma 2.1 }(Bang). {\em For any $q\in
\mathbb N$ and any $x, x+h \in [0,1]$,
$$
B_f(x+h) < \max\left\{ B_f(x), e^{-q}\right\} e^{e|h|A(q)}.
\eqno (2.2)
$$
}

\medskip\par\noindent{\em Proof of the lemma:}
We fix $j$ in the range $0\le j \le q-1$ and find $\xi$ between
$x$ and $x+h$ such that
$$
f^{(j)}(x+h) = \sum_{l=0}^{q-j-1} \frac{f^{(j+l)}(x) h^l}{l!} +
\frac{f^{(q)}(\xi) h^{q-j}}{(q-j)!}\,.
$$
Then
\begin{eqnarray*}
\frac{|f^{(j)}(x+h)|}{e^j M_j} &\le& \sum_{l=0}^{q-j-1}
\frac{|f^{(j+l)}(x)|\, |h|^l}{e^{j} M_{j} l!} +
\frac{|f^{(q)}(\xi)| |h|^{q-j}}{e^j M_j (q-j)!} \\ \\
&=& \sum_{l=0}^{q-j-1} \frac{|f^{(j+l)}(x)|}{e^{j+l} M_{j+l}} \,
\frac{M_{j+l}}{M_j} \, \frac{e^l |h|^l}{l!} + e^{-q} \,
\frac{|f^{(q)}(\xi)|}{M_q} \, \frac{M_q}{M_j} \,
\frac{|h|^{q-j} e^{q-j}}{(q-j)!} \\ \\
&\stackrel{(1.3)}\le& B_f(x) \sum_{l=0}^{q-j-1}
\left(\frac{M_q}{M_{q-1}}\right)^l \, \frac{e^l |h|^l}{l!} +
e^{-q} \left( \frac{M_q}{M_{q-1}} \right)^{q-j} \,
\frac{|h|^{q-j} e^{q-j}}{(q-j)!} \\ \\
&<& \max \left\{B_f(x), e^{-q} \right\} \exp \left( e|h|
\frac{M_q}{M_{q-1}} \right).
\end{eqnarray*}
If $ j\ge q $, the same estimate holds for a trivial reason:
$$
\frac{|f^{(j)}(x+h)|}{e^j M_j} \le e^{-q} < \max\left( B_f(x),
e^{-q} \right) e^{e|h|A(q)},
$$
completing the argument. \hfill $\Box$

\medskip\par\noindent{\bf Corollary 2.3 }{\em Suppose $f\in \DCA$.
If
$$
\max_{[0,1]} B_f \ge e^{-L}, \eqno (2.4)
$$
and
$$
\min_{[0,1]} B_f \le e^{-N}, \eqno (2.5)
$$
then
$$
\sum_{L+1\le j\le N} \frac{M_{j-1}}{M_j} < e\,. \eqno (2.6)
$$
}

\medskip\par\noindent{\em Proof of the corollary:}
Let $B_f(x_N)=e^{-N}$ and $B_f(x_L)=e^{-L}$. By (2.2), the
function $ x\mapsto B_f(x) $ is continuous on $ [0,1] $.
Therefore, we can choose a monotonic sequence $\{x_j\}_{L\le j\le
N} \subset J$ such that $B_f(x_j)=e^{-j}$ for $L\le j\le N$. By
Lemma~2.1,
$$
|x_j-x_{j-1}| > \frac1{e} \frac{M_{j-1}}{M_j},
$$
so that
$$
1 \ge \sum_{L+1\le j\le N} |x_j - x_{j-1}| > \frac{1}{e}
\sum_{L+1\le j \le N} \frac{M_{j-1}}{M_j}\,,
$$
proving the corollary. \hfill $\Box$

\medskip
If the function $f\in\DCA$ has a zero of order at least $N$ at
some point $x_0\in [0,1]$ (that is,
$f(x_0)=f'(x_0)=\,...\,=f^{(N-1)}(x_0)=0$), then $B_f(x_0)\le
e^{-N}$. On the other hand,
$$
\max_{[0,1]} B_f \ge ||f||_{[0,1]}.
$$
Then the corollary says that {\em the order $N$ of any zero of $f$
is bounded from above by the Bang degree $\N$}. This is a
version of a theorem of Borel and Carleman mentioned above. In
particular, the uniqueness part of the Denjoy-Carleman theorem
follows at once: {\em non-trivial functions $f$ from the
quasianalytic Denjoy-Carleman class $\DCA$ cannot have a zero of
infinite order}.

In the non-quasianalytic case, when
$$
\sum_{j=1}^\infty \frac{M_{j-1}}{M_j} <\infty,
$$
rescaling estimate (2.6), we get an upper bound for the function
$f$ near its zeroes of infinite order.

\medskip\par\noindent{\bf Corollary 2.7 }{\em Suppose $f\in \DCA$. Let
$$
f^{(j)}(0) = 0\,, \qquad j\in\mathbb Z.
$$
Then
$$
\sum_{j\ge \log \mathfrak M^{-1}_f(c) +1} \frac{M_{j-1}}{M_j} < ec
$$
where
$$
\mathfrak M_f(c) = \max_{x\in [0,c]}|f(x)|\,.
$$
}

\medskip
Under additional regularity assumptions on the function $A$, Matsaev
and Sodin recently found in \cite{MS} the sharp asymptotics for
$$
 \log \sup_f \mathfrak M_f(c), \qquad c\to 0,
$$
where the supremum is taken over all functions $f\in \mathcal
C_{\{M_j\}}(\mathbb R)$ with $M_0=1$, having the zero of infinite order at
the origin.

\bigskip\centerline{\sc \S 3. Proof of Theorem~A}

\bigskip Now, we consider a sequence of ``norms'' obtained from the
remainders:
$$
b_{f, n}(x) = \max_{j\ge n}\, \frac{|f^{(j)}(x)|}{e^j M_j} =
e^{-n}B_{f^{(n)}}(x),
$$
here $f^{(n)}$ is considered in the class $\mathcal C_{\{M_n,
M_{n+1},\, ... \,\}}([0,1])$. List some properties of this sequence:
\begin{itemize}
\item[(i)] $b_{f, n}(x) \le e^{-n}$;
\item[(ii)] $B_f(x) = b_{f, 0}(x) \ge b_{f, 1}(x) \ge\, ... \,
\ge b_{f, n}(x) \ge \,... $;
\item[(iii)] if $f^{(n)}(x^*)=0$, then $b_{f, n}(x^*)=b_{f, n+1}(x^*)$;
\item[(iv)] the function $b_{f,n}$ satisfies the estimate
$$
b_{f,n}(x+h) \le \max\left\{ b_{f,n}(x), e^{-q-n}\right\} e^{e|h|A(q)},
$$
for every $q\in \mathbb N$ and every $x, x+h\in [0,1]$.
\end{itemize}
From the last property we conclude that
\begin{itemize}
\item[(v)] the function $x\mapsto b_{f,n}(x)$ is continuous on $[0,1]$;
\item[(vi)] $b_{f, n}(x+h)< e^{-j+1}$ provided $b_{f, n}(x)\le
e^{-j}$ and $e|h|A(j)\le 1$.
\end{itemize}
The latter is interesting only for $j>n$.

After these preliminaries we start the proof. Let $x_*\in [0,1]$ be
the maximum point of $B_f$:
$$
\max_{[0,1]} B_f = B_f(x_*).
$$
First, we consider a special case, when $x_*$ is one of the end-points of
$[0,1]$. Without loss of generalities, suppose that $x_*=0$.
Let
$$
0<\xi_1 \le \, ... \, \le \xi_N
$$
be the zeroes of $f$ on $[0,1]$ counted with their multiplicities.
Applying Rolle's theorem, we find another $N$-point set
$\{x_j\}_{0\le j\le N-1}$, such that
$$
f^{(j)}(x_j) = 0, \qquad \hbox{for} \quad 0\le j \le N-1,
$$
$x_0=\xi_1$, and $x_j\le \xi_{j+1}$. Then we can paste together
the functions $b_{f, j}$ with different $j$ and define the new
function
$$
b_{f}(x) = \left\{
\begin{array}{ll}
b_{f, 0}(x), & 0 \le x < x_0 \\
b_{f, 1}(x), & x_0\le x < x_1 \\
... & ... \\
b_{f, N-1}(x), & x_{N-2} \le x < x_{N-1} \\
b_{f, N}(x), & x\ge x_{N-1}\,.
\end{array}
\right.
$$
This is a continuous function with the following properties:
$$
b_f(0) = B_f(0) \ge ||f||_{[0,1]},
$$
$$
b_f(x) \le e^{-N}, \qquad \hbox{for} \quad x\ge x_{N-1},
$$
and $ b_f(x+h) < e^{-j+1}$ provided that $b_f(x)\le e^{-j}$ and
$eh A(j)\le 1$, $h>0$.

Computing as above, in the proof of Corollary~2.3, the
number of level crossings of the function $b_f$ we obtain
$$
1 > \frac1{e} \sum_{j=K_f+1}^N \frac{M_{j-1}}{M_j}
$$
where $K_f= \left\lfloor |\log ||f||_{[0,1]}\right\rfloor$.
Thereby $N\le \N$, completing the proof in the special case.

Now, consider the general case. If $x_*$ is not the end-point of
$[0,1]$, then $x_*$ splits $[0,1]$ into two subintervals $J_1$ and $J_2$
on which $f$ has $N_1$ and $N_2$ zeroes respectively, $N_1+N_2=N$.
By the special case proven above (rescaling the argument of $f$) we
have
$$
|J_l| > \frac1{e} \sum_{j=K_f+1}^{N_l} \frac{M_{j-1}}{M_j},
\qquad l=1,2.
$$
At last, making use of the logarithmic convexity of the sequence
$\{M_j\}$ we obtain
\begin{eqnarray*}
1 &=& |J_1| + |J_2| > \frac1{e} \left\{ \sum_{j=K_f+1}^{N_1}
+\sum_{j=K_f+1}^{N_2} \right\} \frac{M_{j-1}}{M_j} \\ \\
&\ge& \frac1{e} \left\{ \sum_{j=K_f+1}^{N_1} +
\sum_{j=N_1+1}^{N_1+N_2} \right\} \frac{M_{j-1}}{M_j}= \frac1{e}
\sum_{j=K_f+1}^{N} \frac{M_{j-1}}{M_j}
\end{eqnarray*}
whence $N\le \N$. This completes the proof of Theorem A. \hfill $\Box$

\bigskip\centerline{\sc \S 4. Proof of Theorem~B}

\medskip We put $m_0=M_0=1$ and
$$
m_j = \frac{M_j}{j!}, \qquad j\in \mathbb
N.
$$
We shall prove Theorem~B in three steps. First, we prove a
preliminary version of estimate (1.8) with a remainder term:
$$
||f||_I \le \left( \frac{2e |I|}{|E|} \right)^{n} ||f||_E +
m_{n+1}|I|^{n+1} \eqno (4.1)
$$
for each $n\in \mathbb N$. Then we shall show that if the interval
$I$ is sufficiently short, then the remainder $m_{2\N} |I|^{2\N}$
is smaller than the norm $||f||_I$ we are estimating. Combined
with (4.1) this gives us estimate (1.8) for short intervals $I$.
At the last step, we shall extend estimate (1.8) to arbitrary
intervals $I\subset [0,1]$.

\medskip\par\noindent{\bf Proof of estimate (4.1):}
We choose well-spaced points $\{x_j\}_{j=1}^{n+1} \subset E$,
$$
x_1<x_2<...<x_{n+1}, \qquad \min_j (x_{j+1}-x_j) \ge
\frac{|E|}{n},
$$
and set
$$
Q(x)=\prod_{j=1}^{n+1} (x-x_j).
$$
Then for $x\in I$
$$
f(x) = \sum_{j=1}^{n+1} \frac{f(x_j) Q(x)}{Q'(x_j)(x-x_j)} +
\frac{f^{(n+1)}(\xi ) Q(x)}{(n+1)!}\,, \qquad \xi = \xi_x \in I\,.
\eqno (4.2)
$$
This is a well-known version of the Lagrange interpolation
formula. The proof goes as follows: fix $x\in I$ and consider the
function $$G(t) = Q(t)R(x) - Q(x)R(t)$$ where $R(t)$ is the
remainder; i.e. the difference between $f$ and the Lagrange
interpolation polynomial of degree $n$ with the nodes at
$\{x_j\}$. The function $G(t)$ has at least $n+2$ zeroes on $I$:
it vanishes at $n+1$ points: $t=x_j$ and also at $t=x$. Therefore,
the derivative $G^{(n+1)}(t)$ vanishes at least once on $I$:
\begin{eqnarray*}
0=G^{(n+1)}(\xi) &=& Q^{(n+1)}(\xi) R(x) - Q(x) R^{(n+1)}(\xi) \\ \\
&=& (n+1)!R(x) - Q(x) f^{(n+1)}(\xi),
\end{eqnarray*}
proving (4.2).

Then using the estimates $||Q||_I \le |I|^{n+1}$ and
$$
\Big|\Big| \frac{Q(x)}{x-x_j} \Big|\Big|_I \le |I|^{n},
$$
we get
$$
||f||_I \le \left( \sum_{j=1}^{n+1} \frac{1}{|Q'(x_j)|}\right)
|I|^{n} ||f||_E + m_{n+1}|I|^{n+1}\,.
$$

Further,
\begin{eqnarray*}
|Q'(x_j)| &=& (x_j-x_{j-1}) ... (x_j-x_1) (x_{j+1}-x_j) ...
(x_{n+1}-x_j)
\\ \\
&\ge& \frac{(j-1)!(n+1-j)!}{n^n} |E|^{n} >
\frac{(j-1)!(n+1-j)!}{n!} \left(\frac{|E|}{e}\right)^{n}\,,
\end{eqnarray*}
so that
$$
\sum_{j=1}^{n+1} \frac{1}{|Q'(x_j)|} < \left(
\frac{2e}{|E|}\right)^{n}\,,
$$
and (4.1) follows. \hfill $\Box$

\medskip We shall use estimate (4.1) with $n=2\N-1$.

\medskip\par\noindent{\bf Lemma 4.3. }{\em Suppose
$$
m_{2\N} |I|^{2\N} \le e^{-2\N(3+\gamma (2\N)) }
\eqno (4.4)
$$
Then
$$
m_{2\N} |I|^{2\N} < \frac1{2} ||f||_I\,. \eqno (4.5)
$$
}

\medskip Now, combining estimates (4.1) and (4.5) we get

\medskip\par\noindent{\bf Corollary 4.6. }{\em Suppose the interval
$I$ is short; i.e. estimate (4.4) is valid. Then
$$
||f||_I \le \left( \frac{2e |I|}{|E|}\right)^{2\N} ||f||_E\,.
\eqno (4.7)
$$
}

\medskip\par\noindent{\bf Proof of Lemma 4.3} follows from two claims:

\medskip\par\noindent{\bf Claim 4.8 }{\em Estimate (4.4) yields
$$
m_{2\N} |I|^{2\N} \le e^{-4\N} m_k |I|^k \eqno (4.9)
$$
for each $k$, $0\le k\le \N$.}

\medskip\par\noindent{\bf Claim 4.10 }{\em
There exists $k$, $0\le k\le \N$, such that
$$
m_k \left( \frac{|I|}{2}\right)^k \le e^{(2+\frac1{e}) \N} \left[
||f||_I + m_{2\N}\left(\frac{|I|}{2}\right)^{2\N}\right]. \eqno
(4.11)
$$
}

\medskip First, we finish off the proof of Lemma 4.3 and then will
prove the claims. Putting the claims together, we get
\begin{eqnarray*}
m_{2\N} |I|^{2\N} &\le e^{-4\N} \cdot 2^{\N} \cdot e^{(2+\frac1{e}) \N}
\left[ ||f||_I + m_{2\N}\left(\frac{|I|}{2}\right)^{2\N}\right] \\ \\
&\le e^{-0.9\N} \left[ ||f||_I + \frac1{4}m_{2\N}|I|^{2\N}\right]
\end{eqnarray*}
whence
$$
2m_{2\N} |I|^{2\N} < \left( e^{0.9} - \frac1{4}\right) m_{2\N}
|I|^{2\N} \le \left( e^{0.9\N} - \frac1{4}\right) \le  m_{2\N}
|I|^{2\N} ||f||_{I},
$$
proving the lemma. \hfill $\Box$

\medskip\par\noindent{\em Proof of Claim 4.8:} is straightforward.
We have
$$
\frac{m_{2\N}}{m_k} |I|^{2\N-k} \stackrel{(4.4)}\le
e^{-(2\N-k)(3+\gamma (2\N) )}\, \frac{m^{k/2\N}_{2\N}}{m_k}\,.
$$
Therefore, we need to estimate the expression
$$
\left( \frac{m^{k/(2\N)}_{2\N}}{m_k}\right)^{\frac1{2\N-k}} =
\frac{m^{\frac{1}{2\N-k} - \frac1{2\N}}_{2\N}}{m^{\frac1{2\N-k}}_k}
$$
Taking the logarithm and setting $a(s)=\frac{A(s)}{s}$, that is
$m(k)=a(1)a(2)\,...\,a(k)$, we obtain
$$
\frac1{2\N-k} \sum_{j=k+1}^{2\N} \log a(j) - \frac1{2\N}
\sum_{j=1}^{2\N} \log a(j)
$$
$$
= \frac1{2\N-k}\sum_{j=k+1}^{2\N} [\log a(j) - \log a(1)] -
\frac1{2\N} \sum_{j=1}^{2\N} [\log a(j) - \log a(1)]
$$
$$
= \int_1^{2\N} \frac{a'}{a}\big(s\big) \left\{
\frac1{2\N-k}\sum_{j=k+1}^{2\N} \chi_{[1,j]} - \frac1{2\N}
\sum_{j=1}^{2\N} \chi_{[1,j]} \right\}\big(s\big)\, ds \eqno
(4.12)
$$
where $\chi_{[a,b]}(s)$ is the indicator function of the interval
$[a,b]$. Since
$$
0\le \Big\{ \, ... \, \Big\}\big(s\big) < \frac{s}{2\N},
$$
then we get
$$
\hbox{the\ RHS\ of\ }(4.12) < \frac1{2\N}\int_{1}^{2\N}
\frac{sa'(s)}{a(s)} \,ds \le \sup_{s\ge 1}\frac{sa'(s)}{a(s)}
\stackrel{(1.7)}= \gamma (2\N)-1\,.
$$
Therefore
$$
\frac{m^{k/(2\N)}_{2\N}}{m_k} < e^{(2\N-k)(\gamma (2\N)-1)}
$$
and
\begin{eqnarray*}
\frac{m_{2\N}}{m_k} |I|^{2\N-k}
&<& \exp \left\{ -(2\N-k)\left(3+\gamma (2\N) -
\gamma (2\N) +1\right) \right\}
\\ \\
&=& \exp\left\{ -4(2\N-k)\right\} \le \exp \big\{-4\N \big\},
\end{eqnarray*}
proving the claim. \hfill $\Box$

\medskip\par\noindent{\em Proof of Claim 4.10: }
Let $c_I$ be the centre of the interval $I$, and let
$$
P_{2\N-1} (x)= \sum_{j=0}^{2\N-1}
\frac{f^{(j)}(c_I)}{j!}(x-c_I)^j
$$
be the Taylor polynomial of $f$ at $c_I$. Then for $x\in I$
$$
f(x) = P_{2\N-1}(x) + \frac{f^{(2\N)}(\xi)}{(2\N)!}(x-\xi)^{2\N}\,,
\qquad \xi=\xi_x\in I\,,
$$
so that
$$
||P_{2\N-1}||_I \le ||f||_I +
m_{2\N}\left(\frac{|I|}{2}\right)^{2\N}\,.
$$

For an arbitrary polynomial $S$ we have
$$|S^{(k)}(0)|\le ({\rm deg} S)^k ||S||_{[-1,1]}. $$
This is a relatively simple special case of V.~Markov's
inequalities see e.g. \cite[Chapter~VI, Sections 4.II and
4.III]{Mandelbrojt}). Using this inequality, we get
\begin{eqnarray*}
|f^{(k)}(c_I)| &=& |P_{2\N-1}^{(k)}(c_I)| \\ \\
&\le& \left( \frac{2}{|I|}\right)^k (2\N-1)^k ||P_{2\N-1}||_I \\ \\
&\le& \left( \frac{2}{|I|}\right)^k (2\N-1)^k \left[
||f||_I + m_{2\N}\left(\frac{|I|}{2}\right)^{2\N} \right]
\end{eqnarray*}
and
\begin{eqnarray*}
\frac{|f^{(k)}(c_I)|}{e^kM_k} &\le& \frac{1}{m_k} \left(
\frac{2}{|I|} \right)^k \left( \frac{2\N-1}{k} \right)^k
\left[ ||f||_I + m_{2\N}\left(\frac{|I|}{2}\right)^{2\N} \right] \\
&\le& \frac{1}{m_k} \left(\frac{2}{|I|}\right)^k \exp \left(
\frac{2\N-1}{e} \right) \left[ ||f||_I +
m_{2\N}\left( \frac{|I|}{2}\right)^{2\N} \right].
\end{eqnarray*}
Then using Corollary~2.3 from Bang's fundamental lemma and the
definition of the Bang degree $\N$, we obtain that
$$
B_f(c_I) \ge \min_{[0,1]} B_f > e^{-\N-1}\,.
$$
Hence, for {\em at least one} $k$, $0\le k\le \N$,
\begin{eqnarray*}
e^{-\N-1} &<&  \frac{|f^{(k)}(c_I)|}{e^k M_k} \\ \\
&\le&
\frac{1}{m_k} \left( \frac{2}{|I|}\right)^k \exp \left(
\frac{2\N-1}{e}\right) \left[ ||f||_I +
m_{2\N}\left(\frac{|I|}{2}\right)^{2\N}\right],
\end{eqnarray*}
whence
\begin{eqnarray*}
m_k \left(\frac{|I|}{2}\right)^k &\le&
e^{(2+\frac1{e})\N}
\left[ ||f||_I + m_{2\N}\left(\frac{|I|}{2}\right)^{2\N}\right],
\end{eqnarray*}
proving the claim. \hfill $\Box$

\medskip
It remains to spread estimate (4.7) from short to arbitrary
sub-intervals $I\subset [0,1]$. We shall prove a bit more: we show
that if the interval $I\subset [0,1]$ is not short, then
$$
\left( \frac{\Gamma(2\N)|I|}{|E|}\right)^{2\N} ||f||_E \ge 1 \eqno
(4.13)
$$
for any measurable subset $E\subset I$. Since $||f||_I \le
||f||_{[0,1]}\le 1$, this does the job.

We fix a measurable subset  $E\subset I$ where $I$ is not short,
and choose a short sub-interval $I_1\subset I$ such that
$$
|E\cap I_1| \ge |E| \frac{|I_1|}{|I|}, \eqno (4.14)
$$
and
$$
m_{2\N} |I_1|^{2\N} \ge 2^{-2\N}
e^{-2\N\left(3+\gamma (2\N)\right)}.
\eqno (4.15)
$$
Existence of such $I_1$ follows by a straightforward dyadic
argument. Then
\begin{eqnarray*}
1 &\stackrel{(4.15)}\le& \left(2e^{3+\gamma (2\N)}\right)^{2\N}
m_{2\N}|I_1|^{2\N} \\ \\
&\stackrel{(4.5)}\le& \frac{1}{2}
\left(2e^{3+\gamma (2\N)}\right)^{2\N} ||f||_{I_1} \\ \\
&\stackrel{(4.7)}\le& \frac{1}{2} \left(2e^{3+\gamma
(2\N) }\right)^{2\N} \left( \frac{2e |I_1|}{|E\cap I_1|}
\right)^{2\N}
||f||_{E\cap I_1} \\ \\
&\stackrel{(4.14)}<& \left(4e^{4+\gamma (2\N)}
\frac{|I|}{|E|}\right)^{2\N} ||f||_E
\end{eqnarray*}
proving (4.13) and completing the proof of Theorem~B. \hfill
$\Box$

\bigskip\centerline{\sc \S 5. Variations on Bang's theme}

\medskip\par\noindent{\bf 5.1 Bang's differential inequality.}
One can rewrite the fundamental inequality (2.2) as the
differential inequality for the function
$$
L_f(x) = \log \frac1{B_f(x)}.
$$
If $B_f(x)$ is positive, then taking the logarithms in (2.2) and
choosing there $q=[L_f(x)]+1$, we obtain
$$
L_f(x+h)> L_f(x) - e|h|A(L_f(x)+1).
$$
Interchanging $x$ and $x+h$, we arrive at

\medskip\par\noindent{\bf Corollary 5.1.1 }{\em Suppose $f\in \DCA$. Then
the function $L_f: [0,1]\to [0,\infty]$ is continuous and
$$
|L_f(x+h)-L_f(x)| < e|h|A(L_f(x)+1)
$$
whenever $x, x+h\in [0,1]$ and the values $L_f(x)$, $L_f(x+h)$ are
finite.
In particular, if $L'_f(x)$ exists, then
$$
|L_f'(x)| < eA(L_f(x)+1). \eqno (5.1.2)
$$
}

It is remarkable that the function $L_f$ satisfies a simple
differential inequality. Integrating this inequality, we get a
reformulation of Corollary~2.3:

\medskip\par\noindent{\bf Corollary 5.1.3 }{\em Suppose $f\in \DCA$. Then
$$
\int_{L_*+1}^{L^*+1} \frac{ds}{A(s)} <e
$$
where
$$
L_*=\min_{x\in [0,1]} L_f(x)\,, \qquad \hbox{and} \quad
L^*=\max_{x\in [0,1]} L_f(x).
$$
}

\medskip\par\noindent{\bf 5.2 One-sided version of the Denjoy-Carleman 
theorem.}

\medskip\par\noindent{\bf Theorem 5.2.1 }{\em
Suppose $f\in C^\infty([0,1])$ and
$$
\min_{[0,1]} f^{(j)}\ge -M_j \qquad j\in \mathbb Z_+, \eqno
$$
where the sequence $\{M_j\}$ satisfies the quasianalyticity
condition (1.5). If all derivatives of $f$ are non-negative at the
origin, then they are non-negative everywhere on $[0,1]$.}

\medskip Recall that $C^\infty ([0,1))$-functions
with all derivatives positive everywhere on $[0,1)$ are called
{\em absolutely monotonic}. By the classical Bernstein theorem,
every absolutely monotonic function on $[0,1)$ has an analytic
extension to the unit complex disc.

\medskip Under a somewhat stronger assumption $f\in \DCA$,
this result was conjectured by Borel (see \cite[p.74]{Carleman})
and proved by Tacklind \cite{Tacklind} and Bang \cite{Bang}.

For the proof of Theorem~5.2.1, we set
$$
B_f^-(x) = \max_{j\ge 0} \,\frac{\max\{-f^{(j)}(x), 0\}}{e^j M_j}.
$$
Repeating verbatim the proof of Lemma~2.1, we obtain that
$$
B_f^-(x+h) < \max\left\{B_f^- (x), e^{-q} \right\} e^{ehA(q)}
$$
for every $q\in\mathbb N$, every $x,x+h\in [0,1]$, $h>0$. Then
Theorem~5.2.1 readily follows from this estimate. \hfill $\Box$

\medskip\par\noindent{\bf 5.3 Non-extendable quasianalytic
functions. } If $f$ is a real analytic function on a closed
interval $J$ (that is $f\in\mathcal C_{\{K^jj!\}}(J)$ with some
constant $K$), then $f$ always has a real analytic extension on a
larger interval $J'\supset J$. In contrast, the Tacklind-Bang
theorem combined with the Bernstein theorem give us examples of
quasianalytically smooth functions defined on a closed interval
which do not have a quasianalytically smooth extension; i.e. a
smooth extension which belongs to a (probably, different)
Denjoy-Carleman quasianalytic class on a larger interval.

We fix a logarithmically convex sequence $\{M_j\}$ satisfying
(1.5) and such that
$$
\lim_{j\to \infty} \left(\frac{M_j}{j!}\right)^{1/j} = +\infty,
$$
and choose a positive sequence $\{c_j\}$ such that
$$
\lim_{j\to \infty} c_j^{1/j} =1,
$$
and
$$
\sum_{j=1}^\infty j^n c_j \le M_n
$$
for all $n\in \mathbb Z_+$. For example, if $M_n=n!(\log n)^n$,
then one can take
$$c_j=\exp\left[- \frac{Cj}{\log(j+e)}\right]$$ with a proper choice
of a positive constant $C$.

Then consider an even function
$$
f(x)= \sum_{k=0}^\infty c_{2k} x^{2k}
$$
which is is analytic in $(-1,1)$ and belongs to the quasianalytic
class $\mathcal C_{\{M_j\} }$ on the segment $[-1,1]$. This
function has no quasianalytically smooth extension on a larger
interval. Otherwise, the extension would be an even function (by
the Denjoy-Carleman theorem), and by the theorems of Tacklind-Bang
and Bernstein it would have an analytic extension to a disk of
radius larger than one. Clearly, this is impossible since the
radius of convergence of the Taylor series which represents $f$
equals one.

This construction answers the question raised by P.~Milman
\footnote{A. Borichev indicated another construction of a
non-extendable quasianalytic function. He considers the absolutely
convergent series
$$
f(z) = \sum_n \frac{\delta_n}{z-\lambda_n}, \qquad \lambda_n =
1+\epsilon_n-i\delta_n,
$$
with
$$
\epsilon_n\downarrow 0, \qquad \epsilon_{n+1}/\epsilon_n\downarrow
0, \qquad \delta_n\downarrow 0, \qquad \delta_{n+1}/\delta_n
\downarrow 0, \qquad \delta_n/\epsilon_n \downarrow 0
$$
and using some results from \cite{Beurling} shows that under a
special choice of these sequences the function has no
quasianalytic extension to any larger interval $ [-1-\gamma,
1+\gamma] $. }.

\medskip\par\noindent{\bf 5.4 Bang's original version of the
fundamental inequality. } Mention that Bang proved his results
without assumption of the logarithmic convexity of the sequence
$\{M_j\}$. He assumed that $M_j$ is the upper bound for the
$|f^{(j)}(x)|$ on the closed interval $J$ and the sequence
$\{M_j\}$ increases so rapidly that $M_j^{1/j}\to \infty$. Then
there exists a unique largest logarithmically convex minorant
$M^c_j\le M_j$. The equation $M_j^c=M_j$ is satisfied for
infinitely many integers $j$, in particular for $j=0$ (see e.g.
\cite[Chapter~1]{Mandelbrojt}). This set of integers is denoted by
$\mathbb P$. Then Bang defines the ``norm''
$$
B_f(x) = \inf_{p\in \mathbb P}\, \max \left\{e^{-p}, \max_{0\le
j\le p} \frac{|f^{(j)}(x)|}{e^jM^c_j}\right\}
$$
and proves that if $B_f(x)\ge e^{-q}$, $q\in \mathbb N$, then
$$
B_f(x+h) < B_f(x) e^{e|h|A^c(q)}
$$
where
$$
A^c(q) = \frac{M^c_{q}}{M^c_{q-1}}.
$$
From here he deduces a more general version of Corollary~2.3 which
already contains the result of Borel and Carleman formulated in
the Introduction.

\medskip\par\noindent{\bf 5.5 Propagation of smallness for
quasianalytically smooth functions.} Here, we give a simple
``global corollary'' to Theorem~B excluding the degree $\N$ from
estimate (1.8). We assume that $\DCA$ is a regular quasianalytic
Denjoy-Carleman class of functions, that is $A: [1,\infty)\to
(0,\infty)$ is a non-decreasing $C^1$-function such that the
integral $ \int^\infty A^{-1}(s)\,ds $
is divergent and
$$
\gamma = \sup_{s\ge 1} \frac{sA'(s)}{A(s)} < \infty.
$$

We set
$$
\Omega (t) \stackrel{\rm def}= \exp\left[- \frac1{e}\int_1^{\log(e/t)}
\frac{ds}{A(s)} \right]\,, \qquad 0\le t \le 1\,.
$$
The function $\Omega$ steadily increases on the interval
$[0,1]$, $\Omega (0)=0$, and $\Omega (1)=1$.

A relative smallness of the set $E\subset [0,1]$ will be measured
in the logarithmic scale by the quantity
$$
\alpha (E) = \frac{1}{3}\, \log^{-1} \left( \frac{\Gamma}{|E|}
\right)
$$
where as above $\Gamma=4e^{4+\frac2{e}\gamma }$.

\medskip\par\noindent{\bf Corollary 5.5.1 }{\em
Suppose $\DCA$ is a regular Denjoy-Carleman quasianalytic class,
and suppose that $f\in\DCA$. Then
$$
\Omega \left( ||f||_{[0,1]}\right) \le e \Omega
\left( ||f||_E^{\alpha (E)} \right). \eqno
(5.5.2)
$$
}

\medskip In the real analytic case when $A(s) = Cs$, $C$ is a positive
constant, we have
$$
\Omega (t) = \exp\left[ - \frac1{eC} \log\log \frac{e}{t} \right] =
\left(\log \frac{e}{t} \right)^{-1/(eC)}.
$$
Suppose that $||f||_E \le \epsilon$. Then estimate (5.5.2) gives
us
$$
||f||_{[0,1]} \le e \epsilon^{\alpha (E)
e^{-eC}}. \eqno (5.5.3)
$$
Certainly, estimate (5.5.3) can be obtained by classical means
using a complex extension with control over the uniform norm and
the two-constants-theorem \cite{Lavrentiev}, \cite{Vessella}, or
by an elementary real variable technique \cite{Nadirashvili1}.

However, already in the logarithmic Denjoy-Carleman class when
$A(s)=Cs\log (s+e)$, the Corollary gives a new result.

\medskip\par\noindent{\em Proof of Corollary 5.5.1: }
We have
$$
||f||_{[0,1]} \le \left( \frac{\Gamma}{|E|}\right)^{2\N} ||f||_E
$$
or
\begin{eqnarray*}
1 &\le& \exp\left\{ 2\N \log \left( \frac{\Gamma}{|E|}\right) +
\log \frac1{||f||_{[0,1]}}
\right\} ||f||_E \\ \\
&<& \exp\left\{ 3\N \log\left( \frac{\Gamma}{|E|} \right) \right\}
||f||_E  = \exp\left\{ \frac{\N}{\alpha(E)}\right\} ||f||_E \,,
\end{eqnarray*}
that is
$$
e^{-\N} \le ||f||_E^{\alpha (E)}\,.
$$
Since
\begin{eqnarray*}
\log \Omega \left(||f||_{[0,1]}\right) - \log \Omega (e^{-\N}) &=&
\frac1{e}\int_{\log ||f||^{-1}_{[0,1]}+1}^{\N+1} \frac{ds}{A(s)}
\\ \\
&\le& \frac1{e} \sum_{\log ||f||^{-1}_{[0,1]}<j\le \N}
\frac{M_{j-1}}{M_j} < 1\,,
\end{eqnarray*}
we finally get
$$
\Omega \left( ||f||_{[0,1]} \right) < e \Omega
(e^{-\N}) < e  \Omega \left( ||f||_E^{\alpha (E)}\right),
$$
completing the proof. \hfill $\Box$

\medskip
\par\noindent (F.N.:) Department of Mathematics,

\par\noindent Michigan State
University,

\par\noindent East Lansing, MI 48824,

\par\noindent U.S.A.

\par\noindent {\textit{\small E-mail}: \texttt{\small
fedja@math.msu.edu }}

\medskip
\par\noindent (M.S.:) School of Mathematical Sciences,

\par\noindent Tel Aviv University,

\par\noindent Ramat Aviv, 69978,

\par\noindent Israel

\par\noindent {\textit{\small E-mail}: \texttt{\small
sodin@post.tau.ac.il }}

\medskip
\par\noindent (A.V.:) Department of Mathematics,

\par\noindent Michigan State University,

\par\noindent East Lansing, MI 48824,

\par\noindent U.S.A.

\par\noindent {\textit{\small E-mail}: \texttt{\small
volberg@math.msu.edu }}

\end{document}